\documentclass{amsart}
\usepackage{amsmath, amsthm, amssymb}
\usepackage{verbatim,url,xcolor}
\usepackage{url}
\usepackage{xcolor}

\renewcommand{\d}{\partial}
\newcommand{\dbar}{\overline{\partial}}
\newcommand{\ddbar}{\sqrt{-1}\d\overline{\d}}

\newtheorem{thm}{Theorem}
\newtheorem{prop}[thm]{Proposition}

\newtheorem{ex}[thm]{Example}

\theoremstyle{definition}

\newtheorem{rem}[thm]{Remark}

\renewcommand{\[}{\begin{equation}}
\renewcommand{\]}{\end{equation}}

\newcommand{\vep}{\varepsilon}

\newcommand{\Rm}{\mathrm{Rm}}

\newcommand{\RR}{\mathbb{R}}
\newcommand{\CC}{\mathbb{C}}

\newcommand{\PP}{\mathbb{C}P}

\newcommand{\Rc}{\mathrm{Ric}}

\newcommand{\BK}{\mathrm{BK}}
\newcommand{\diam}{\mathrm{diam}}
\newcommand{\sn}{\mathrm{sn}}

\title[Diameter rigidity]{Diameter rigidity for K\"ahler manifolds with positive bisectional curvature}
\author[V. V. Datar]{Ved Datar}
\author[H. Seshadri]{Harish Seshadri}
\address{Department of Mathematics, Indian Institute of Science, Bangalore, India - 560012}
\email{vvdatar@iisc.ac.in, harish@iisc.ac.in}

\thanks{Research supported in part by the Infosys Young Investigator award.}

\begin{document}
\begin{abstract}
We prove that a K\"ahler manifold with positive bisectional curvature and maximal diameter is isometric to the complex projective space with the Fubini-Study metric.
 \end{abstract}
\maketitle

\section{Introduction}

Let $(M,\omega)$ be a compact K\"ahler manifold of dimension $n$. The {\it bisectional curvature} of $\omega$ along real unit tangent vectors $X,Y$ is defined to be $$\BK(X,Y) = \Rm(X,JX,JY,Y),$$  
where ${\rm Rm}$ denotes the Riemann curvature tensor of the Riemannian metric associated to $\omega$. 
In this note we will be concerned with K\"ahler manifolds $(M,\omega)$ satisfying 
\begin{align}\label{eq:bk}
\BK \geq  1,
\end{align} 
i.e., $\BK(X,Y) \ge 1$ for all real unit tangent vectors $X,Y$.  

%
%A Laplacian comparison theorem for K\"ahler manifolds with  bisectional curvature bounded below was proved in \cite{LW} and, as a corollary, the following upper bound  for the diameter was established:
\begin{thm}[Li-Wang \cite{LW}]\label{LW}
If $(M,\omega)$ is a compact K\"ahler manifold satisfying $\BK \ge 1$, then
$$\diam(M) \ \leq \ \diam(\PP^n,\omega_{\PP^n})= \frac{\pi}{\sqrt{2}} .$$
\end{thm} 

We address the case of equality in Theorem \ref{LW} in this note:    
\begin{thm}\label{main}
Let $(M^n,\omega)$ be a compact K\"ahler manifold satisfying 
$\BK \ge 1.$ 
If $$\diam(M,\omega) = \diam(\PP^n,\omega_{\PP^n}) ,$$ then $(M,\omega)$ is isometric to $(\PP^n,\omega_{\PP^n}).$
\end{thm}
\vspace{2mm}
Here the Fubini-Study metric is normalized so that
  $$\int_{\PP^n}\omega_{\PP^n}^n = (2\pi)^n, \ {\rm equivalently} \ {\rm Ric} =(n+1) \omega_{\PP^n}.$$
\begin{rem} In \cite{LW}, the  diameter bound is stated to be $\pi/2$. This is due to a
 different normalization for the Hermitian extension of the Riemannian metric.  
\end{rem}

The diameter bound in Theorem \ref{LW} is analogous to the classical Bonnet-Myers diameter bound for compact Riemannian manifolds with positive Ricci curvature. However, one cannot relax the curvature assumption to a positive Ricci lower bound  in the K\"ahler case: consider ${\mathbb C}P^1$ endowed with the round metric of curvature $\frac{1}{n+1}$ and take the product metric on the $n$-fold product
$$M= {\mathbb C}P^1 \times ... \times {\mathbb C}P^1.$$
The Ricci curvature of $M$ satisfies $\Rc = (n+1)\omega$, but $${\rm diam}(M) =\sqrt{\frac{n}{n+1}} \pi >
\frac{\pi}{\sqrt{2}}$$ if $n\ge2$. \vspace{2mm}
 
In the Riemannian case, the equality case of the Bonnet-Myers diameter bound is addressed by the well-known maximal diameter theorem of Cheng. Our main result can be regarded as the  K\"ahler analogue of  Cheng's theorem.

%Note that a volume rigidity result was recently proved in \cite{Z1} under only a Ricci curvature lower bound. 
Theorem \ref{main} has been established under additional assumptions in  \cite{GY} and \cite{TY}.  In \cite{GY}, the authors construct a totally geodesic $\PP^1$ with sectional curvature $2$ and use this to  show  that rigidity holds if $\int_M\omega^n > \pi^n$.  In \cite{TY}, the authors assume that  there are complex submanifolds $P$ and $Q$ of $M$ with ${\rm dim}(P) +{\rm dim}(Q) =n-1$ and  prove rigidity. An eigenvalue comparison theorem is the main ingredient in their proof.

Our strategy for proving Theorem \ref{main} is to establish a monotonicity formula for a function arising from Lelong numbers of positive currents on $\PP^n$. In \cite{Lott}, the $\partial \bar \partial$-comparison theorem of \cite{TY} is reformulated as asserting the positivity of a certain $(1,1)$-current and this is the current we work with. 
\section*{Acknowledgements}
We would like to thank Vamsi Pingali for his interest in the work and helpful discussions. We would also like to thank John Lott for useful comments on the first draft of the paper. 
% in \cite{DSS}.  An interesting question that we wish to address in a sequel is if the corresponding almost diameter rigidity holds. 

\section{Lelong numbers and a monotonicity formula on $\PP^n$} 
\noindent Let $M$ be a K\"ahler manifold. In what follows, we frequently use the real operator  $$~d^c = \frac{\sqrt{-1}}{2\pi}(\dbar-\d).$$
Note that $$dd^c = \frac{1}{\pi}\ddbar.$$

\noindent If $T$ is a non-negative current on a $M$ such that $$T  =  dd^c\varphi$$ in a neighbourhood of a point $q\in M$, then the {\it Lelong number} of $T$ at $q$ is defined as $$\nu(T,q) :=\lim_{r\rightarrow 0^+}\frac{\sup_{B_{\CC^n}(0,r)}\varphi(z)}{\log r},$$ where $z$ is a holomorphic coordinate in a neighbourhood of $q$ such that $z(q) = 0$. It is not difficult to see (for instance using the maximum principle) that the quotient on the right is increasing in $r$, and hence the limit $\nu(T,q)$ exists and is moreover non-negative and independent of the choice of holomorphic coordinates. Note  that the normalization is chosen so that if $V$ is a smooth hypersurface with defining function $f$, and $[V]$ denotes the current of integration along $V$, then by the Poincare-Lelong equation, $[V] = dd^c\log|f|$, and so $\nu([V],q) = 1$ for any point $q\in V$.  

The following proposition is well known (cf. \cite[pg.\ 164-165]{Dem-book}), but since the proof of our main theorem has a precise dependence on the constants involved, we provide a proof for the convenience of the reader. 

\begin{prop} \label{prop:lelong-limit}Suppose $T = dd^c\varphi$ as above in a neighbourhood of $q$ with holomorphic coordinates $z=(z^1,\cdots,z^n)$ such that $z(q) = 0$. We then have
$$\nu(T,q) = \lim_{r\rightarrow 0^+}\frac{1}{\pi^{n-1}r^{2n-2}}\int_{B_{\CC^n}(0,r)}T\wedge \omega_{\CC^n}^{n-1},$$ where $B_{\CC^n}(0,r)$ is the ball of radius $r$ around the origin with respect to the Euclidean metric $\omega_{\CC^n} = \frac{\sqrt{-1}}{2}\partial\bar\partial |z|^2$. 
\end{prop}Note that quantity on the right above  is increasing in $r$ (cf.\ \cite[pg.\ 390]{GH}), and hence the limit in particular exists.

\begin{proof} First suppose that $\varphi$ is smooth. We let 

$$\nu(dd^c\varphi,0,t) := \frac{1}{\pi^{n-1}t^{2n-2}}\int_{B_{\CC^n}(0,t)}dd^c\varphi\wedge \omega_{\CC^n}^{n-1},$$ 

$$\mu_t(\varphi):= \frac{1}{\sigma_{2n-1}}\int_{\mathbb{S}^{2n-1}}\varphi(t,\theta)\,d\sigma(\theta),$$ 

where $\sigma_{2n-1} = 2\pi^n/(n-1)!$ is the volume of the unit sphere in $\mathbb{S}^{2n-1}\subset \CC^n$, and  $d\sigma$ is the standard Riemannian measure on $\mathbb{S}^{2n-1}$
% so that $\mu_t(\varphi)$ is the average of $\varphi$ on $|z|=t$. 
Let $\mathbb{S}^{2n-1}_t$  be the sphere of radius $t$ centred at the origin, $d\sigma_t$ the Riemannian measure on it and let $\partial\varphi/\partial\nu$ be the normal derivative of $\varphi$. Differentiating in $t$, 
\begin{align*}
\frac{d\mu_t(\varphi)}{dt} &=\frac{1}{\sigma_{2n-1}}\int_{\mathbb{S}^{2n-1}}\frac{\partial \varphi}{\partial t}(t,\theta)\,d\sigma\\
&=\frac{1}{\sigma_{2n-1}t^{2n-1}}\int_{\mathbb{S}_t^{2n-1}} \frac{\partial\varphi}{\partial\nu}\,d\sigma_t\\
&=\frac{2}{\sigma_{2n-1}t^{2n-1}}\int_{B_{\CC^n}(0,t)}\Delta_{\dbar} \varphi \frac{\omega_{\CC^n}^n}{n!}\\
&=\frac{2}{\sigma_{2n-1}t^{2n-1}}\int_{B_{\CC^n}(0,t)}\ddbar \varphi \wedge \frac{\omega_{\CC^n}^{n-1}}{(n-1)!}\\
&=\frac{2\pi}{\sigma_{2n-1}(n-1)!}\cdot\frac{1}{t^{2n-1}}\int_{B_{\CC^n}(0,t)}dd^c\varphi \wedge \omega_{\CC^n}^{n-1}\\
&= \frac{\nu(T,0,t)}{t}.
\end{align*}

Note that in the third line we have the $\dbar$-Laplacian $\Delta_{\dbar}$, and hence the factor of $2$ on application of Green's formula. Integrating the above equality from $r$ to $1$, we obtain the so-called Jensen-Lelong formula (cf. \cite[pg.\ 163]{Dem-book}): 

$$\mu_1(\varphi) - \mu_r(\varphi) = \int_{r}^1\nu(dd^c\varphi,0,t)\,\frac{dt}{t}.$$ 

By regularization, the above equality also holds for a general, possibly non-smooth, plurisubharmonic function $\varphi$. Changing variables $s = \log t$ and dividing by $\log r$ we have $$\frac{\mu_r(\varphi)}{\log r} = \frac{\mu_1(\varphi)}{\log r} - \frac{1}{\log r}\int_{\log r}^0\nu(dd^c\varphi,0,e^s)\,ds,$$ and letting $r\rightarrow 0^+$ we obtain $$\lim_{r\rightarrow 0^+}\nu(T,0,r) = \lim_{r\rightarrow 0^+}\frac{\mu_r(\varphi)}{\log r}.$$ Next proceeding as in \cite[pg.\ 165]{Dem-book}, by Harnack inequality and maximum principle, we have that $$ \lim_{r\rightarrow 0^+}\frac{\mu_r(\varphi)}{\log r} = \lim_{r\rightarrow 0^+}\frac{\sup_{z\in \partial B_{\CC^n}(0,r)}\varphi(z)}{\log r} =  \lim_{r\rightarrow 0^+}\frac{\sup_{z\in B_{\CC^n}(0,r)}\varphi(z)}{\log r}.$$
\end{proof}We require the following modification, which as far as we can tell, seems to be new.

\begin{prop}\label{prop:mono}.
Let $T$ be a non-negative current on $\PP^n$ in a K\"ahler class, and $q\in \PP^n$. Then $$\Theta(T,q,r):= \frac{1}{(2\pi)^{n-1}\sin^{2n-2}(r/\sqrt{2})}\int_{B_{\PP^n}(q,r)}T\wedge \omega_{\PP^n}^{n-1}$$ is increasing in $r$. Here $B_{\PP^n}(q,r)$ is the ball of radius $r$ with respect to $\omega_{\PP^n}$. Moreover, we also have that 
\begin{equation}\label{eq:gen-lelong-limit}
\lim_{r\rightarrow 0^+}\Theta(T,q,r) = \nu(T,q).
\end{equation} \end{prop}  

Note that the factor in the denominator is precisely the volume of a ball of radius $r$ in $\PP^{n-1}$ with respect to the Fubini-Study metric $\omega_{\PP^{n-1}}$ upto a factor of $(n-1)!$.
 \begin{proof}
 Let us first assume that $T$ is a smooth $(1,1)$ K\"ahler form. We use homogenous coordinates $[\xi_0:\xi_1:\cdots:\xi_n]$ on $\PP^n$ with $q = [1:0:\cdots:0]$, and the usual in-homogenous coordinates $Z_i = \frac{\xi_1}{\xi_0}$ on $\xi_0\neq 0$. Then $$\omega = \ddbar \log|\xi|^2 = \ddbar\log(1+|Z|^2).$$ We then compute 
\begin{align*}
\Theta(T,q,r)  &= \frac{1}{2^{n-1}\sin^{2n-2}{(r/\sqrt{2}})}\int_{B_{\PP^n}(q,r)}T\wedge (dd^c\log|\xi|^2)^{n-1}\\
&=\frac{1}{2^{n-1}\sin^{2n-2}(r/\sqrt{2})}\int_{\partial B_{\PP^n}(q,r)} T\wedge d^c\log(1+|Z|^2)\wedge (dd^c\log(1+|Z|^2))^{n-2}.
\end{align*}
Now, it is well known fact that $$\cos^2\frac{d_{\PP^n}(q,Z)}{\sqrt{2}} = \frac{|\xi_0|^2}{|\xi|^2} = \frac{1}{1+|Z|^2}.$$ For instance exploiting the $U(n)$ symmetry one needs to check this only for $\PP^1$ which can be done easily. We then have that for any $Z\in \partial B_{\PP^n}(q,r)$, $$d^c\log(1+|Z|^2) = \frac{|Z|^2}{1+|Z|^2}d^c\log|Z|^2 = \sin^2\Big(\frac{r}{\sqrt{2}}\Big)d^c\log|Z|^2.$$ Putting this back in the formula above we have that 
\begin{equation}\label{eq:altformula}
\Theta(T,q,r)  = \frac{1}{2^{n-1}}\int_{\partial B_{\PP^n}(q,r)}T\wedge d^c\log|Z|^2\wedge (dd^c\log|Z|^2)^{n-2}.
\end{equation} So if $r_1<r_2$, then integrating by parts we have $$\Theta(T,q,r_2) - \Theta(T,q,r_1) = \frac{1}{2^{n-1}}\int_{A_{\PP^n}(q,r_1,r_2)}T\wedge (dd^c\log|Z|^2)^{n-1},$$ where  $A_{\PP^n}(q,r_1,r_2) = B_{\PP^n}(q,r_2)\setminus\overline{B_{\PP^n}(q,r_1)}$. Now if $\mu:\PP^n\dashrightarrow \PP^{n-1}$ is the projection from $q$ to $[\xi_0 = 0]$, then we have $$\Theta(T,q,r_2) - \Theta(T,q,r_1) = \frac{1}{(2\pi)^{n-1}}\int_{A_{\PP^n}(q,r_1,r_2)}T\wedge (\mu^*\omega_{\PP^{n-1}})^{n-1} \geq 0.$$ This proves the monotonicity for smooth currents. For a general positive current $T$ we can proceed by regularization. In fact in our case we can first let $r_1<r_2<R<\pi/\sqrt{2}$. Then $B(q,R)$ is contained in Euclidean ball (of radius $\tan R$) with respect to the in-homogenous coordinates. We can then use the standard convolution to find sequence of smooth non-negative forms $T_j$ converging weakly to $T$. Then since $r_1<r_2<R$, 
\begin{align*}
\Theta(T,q,r_2) - \Theta(T,q,r_1) &=\lim_{j\rightarrow \infty}\Big(\Theta(T_{j},q,r_2) - \Theta(T_{j},q,r_1)\Big) \geq 0.
\end{align*} 
If $r_2 = \pi/\sqrt{2}$, then the result follows by the monotonic convergence. 

Next, to compute the limit, we again first work with smooth Kahler forms. If $T$ is smooth then in formula \eqref{eq:altformula}, we observe that $$d^c\log|Z|^2 = \frac{d^c|Z|^2}{|Z|^2} = \frac{d^c|Z|^2}{\tan^2(r/\sqrt{2})},$$ where notice that $d(q,Z) = r$ implies that $$|Z|^2 = \tan^2\Big(\frac{r}{\sqrt{2}}\Big).$$
Then we have 
\begin{align*}
\Theta(T,q,r)  &= \frac{1}{2^{n-1}}\int_{\partial B_{\PP^n}(q,r)}T\wedge d^c\log|Z|^2\wedge (dd^c\log|Z|^2)^{n-2} \\
&=\frac{1}{2^{n-1}\tan^{2n-2}(r/\sqrt{2})}\int_{B_{\PP^n}(q,r)}T\wedge d^c|Z|^2\wedge  (dd^c|Z|^2)^{n-2}\\
&=\frac{1}{2^{n-1}\tan^{2n-2}(r/\sqrt{2})}\int_{B_{\PP^n}(q,r)}T\wedge (dd^c|Z|^2)^{n-1}\\
&=\frac{1}{\pi^{n-1}t^{2n-2}}\int_{B_{\CC^n}(0,t)} T\wedge \omega_{\CC^n}^{n-1},
\end{align*}
where we integrated by parts in the third line and set $t = \tan(r/\sqrt{2})$, and noted that in terms of the $Z$-coordinates $B_{\PP^n}(q,r) = B_{\CC^n}(0,t)$. Once again by regularization, as above, the above formula holds for general possibly non-smooth currents. Letting $t\rightarrow 0^+$ and applying Proposition \ref{prop:lelong-limit} we obtain \eqref{eq:gen-lelong-limit}. 

%$$\lim_{r\rightarrow 0^+}\frac{1}{\pi^{n-1}r^{2n-2}}\int_{B_{\PP^n}(q,r)}dd^c\varphi\wedge \omega_{\PP^n}^{n-1} = 2\nu(T,q),$$ %and so if $T = \ddbar\varphi$ locally we have $$\lim_{r\rightarrow 0^+}\Theta(T,q,r) = \pi\lim_{r\rightarrow 0^+}\frac{(r/\sqrt{2})^{2(n-1)}}{\sin^{2n-2}(r/\sqrt{2})}\lim_{r\rightarrow 0^+}\frac{1}{\pi^{n-1}r^{2n-2}}\int_{B_{\PP^n}(q,r)}dd^c\varphi\wedge \omega_{\PP^n}^{n-1} = 2\pi\nu(T,q). $$ 

 \end{proof}
 
 \begin{ex}[The ``model" case]\label{rem:model}
{\rm On $\PP^n$ consider the current $T = \ddbar\log|\xi_n|^2 = 2\pi[\xi_n=0]$, and $q = [1:0:\cdots:0]$. We regard this  as the model case for reasons given in Section 3. Then for any $r>0$, 
 \begin{align*}
 \int_{B_{\PP^n}(q,r)} T\wedge \omega_{\PP^n}^{n-1} &= 2\pi \int_{B_{\PP^n}(q,r)\cap\{\xi_n = 0\}}\omega_{\PP^n}^{n-1}\\
 &=2\pi\int_{B_{\PP^{n-1}(q,r)}}\omega_{\PP^{n-1}}^{n-1}\\
 &= (2\pi)^{n} \sin^{2n-2}\Big(\frac{r}{\sqrt{2}}\Big),
 \end{align*}
 and so $\Theta(T,q,r) = 2\pi$ and is independent of $r$. Note that if we consider a modified $$\tilde\Theta(T,q,r) := \frac{1}{(2\pi)^{n-1}r^{2n-2}}\int_{B_{\PP^n}(q,r)}T\wedge \omega_{\PP^n}^{n-1},$$ where we have $r^{2n-2}$ in the denominator as in the usual Euclidean case, then  for $T$ and $q$ as above we would have that $$\tilde \Theta(T,q,r) = 2\pi\frac{\sin^{2n-2}(r/\sqrt{2})}{r^{2n-2}}.$$ It is easy to see that this function is {\em decreasing} in $r$ and not increasing.}
 \end{ex}

\section{Proof of the Theorem}

In \cite{Lott}, Lott introduces the following current: $$T_{\omega,p}:= \omega+\ddbar\psi_p,~\psi_p := \log\cos^2\Big(\frac{d_p}{\sqrt{2}}\Big),$$ where $p$ is some fixed point in $M$ and $d_p$ is the distance function from $p$. Note that {\em a priori}, $T_{\omega,p}$ is only defined (and also smooth) away from the cut-locus of $p$. If $\omega = \omega_{\PP^n}$, and $p = [0:0:\cdots:1]$, then as observed before $$\cos^2\Big(\frac{d_{\omega_{\PP^n},p}}{\sqrt{2}}\Big) = \frac{|\xi_n|^2}{|\xi|^2},$$ and so $$T_{\omega_{\PP^n},p} = \ddbar \log|\xi_n|^2 $$ is precisely the current considered in Remark \ref{rem:model} above. Using the Hessian comparison theorem in \cite{TY}, which {\em apriori} holds only away from the cut-locus, Lott observed that $T$ is in fact a global non-negative current if $\omega$ satisfies \eqref{eq:bk}. 

\begin{proof}[Proof of theorem] Firstly note that by the proof of the Frankel conjecture (cf.\ \cite{SY} and \cite{M}), $M$ is bi-holomorphic to $\PP^n$. So from now on we set $M = \PP^n$. Let $p,q\in \PP^n$ such that $d_{\omega,p}(q) = \pi/\sqrt{2}$. We claim that $\nu(T_{\omega,p},q) = \nu(\omega+\pi dd^c\psi_{\omega,p})\geq 2\pi$. Using normal coordinates, it is in fact enough to show that $$\lim_{\vep\rightarrow 0^+}\frac{\sup_{B(q,\vep)}\psi_{\omega,p}}{\log\vep} \geq 2,$$ since $\omega$ being smooth does not contribute to the Lelong number.

It is more convenient to work with $$\delta_p = \frac{\pi}{2} - \frac{d_p}{\sqrt{2}}.$$ Then $\psi_p = 2\log\sin\delta_p$. Note that by the diameter upper bound we have $\delta_p(z)\geq 0$ for all $z$, and that $\delta_p$ is Lipshitz with constant $1/\sqrt{2}$. Then for any $x\in \PP^n$, $$\delta_p(x) =\leq  \frac{1}{\sqrt{2}}d(q,x),$$ and so $\sup_{B_\CC^n(q,\vep)}\psi_{\omega,p}\leq C + 2\log\vep.$ But then $$\frac{\sup_{B(q,\vep)}\psi_{\omega,p}}{\log\vep} \geq \frac{C}{\log\vep} + 2\xrightarrow{\vep\rightarrow 0^+}2.$$%In particular, $$\psi_p(z)= 2\log \sin\delta_p(z) \leq 2\log\delta_p(z) \leq -\log 2 + 2\log d(q,z),$$ and so $\nu(T_{\omega,p},q) \geq 2\pi.$  The factor of $\pi$ comes about since the singular part of $T$ comes from $\ddbar\psi_p$ which equals $\pi dd^c\psi_p$. 
But then by monotonicity, if $\omega \in c[\omega_{\PP^n}]$, putting $R = \pi/\sqrt{2}$, we have
$$
2\pi c =  \frac{1}{(2\pi)^{n-1}}\int_{\PP^n}T\wedge \omega_{\PP^n}^{n-1} = \Theta(T_{\omega,p},q,R) \geq \lim_{r\rightarrow 0^+} \Theta(T_{\omega,p},q,r) = \nu(T_{\omega,p},q) \geq 2\pi,
$$ and so $c\geq 1$. On the other hand note that the bisectional curvature lower bound gives $$\Rc(\omega) \geq (n+1)\omega,$$ and so $c\leq 1$ since $[\Rc(\omega)] = (n+1)[\omega_{\PP^n}]$, and hence $c = 1$. But then the lower bound on the Ricci curvature, and the $\ddbar$-lemma imply that $\omega$ must be K\"ahler-Einstein and hence isometric to $\omega_{\PP^n}$.
\end{proof}

\end{document}